\documentclass[12pt, reqno]{amsart}
\usepackage{amsfonts}
\usepackage{amssymb}
\usepackage{amsmath}

\newtheorem{theorem}{Theorem}
\theoremstyle{plain}

\newtheorem{corollary}[theorem]{Corollary}

\newtheorem{proposition}[theorem]{Proposition}

\begin{document}
\title[] {$q$-Bernstein polynomials associated with $q$-Stirling numbers and Carlitz's $q$-Bernoulli numbers}
\author{T. Kim}
\address{Taekyun Kim. Division of General Education-Mathematics \\
Kwangwoon University, Seoul 139-701, Republic of Korea  \\}
\email{tkkim@kw.ac.kr}
\author{J. Choi}
\address{Jongsung Choi. Division of General Education-Mathematics \\
Kwangwoon University, Seoul 139-701, Republic of Korea  \\}
\email{jeschoi@kw.ac.kr}
\author{Y.H. Kim}
\address{Young-Hee Kim. Division of General
Education-Mathematics\\
Kwangwoon University, Seoul 139-701, Republic of Korea  \\}
\email{yhkim@kw.ac.kr}
%\thanks{$*$ correponding author}
\maketitle

{\footnotesize {\bf Abstract} \hspace{1mm} {Recently, T. Kim([4])
introduced $q$-Bernstein polynomials which are different
$q$-Bernstein polynomials of Phillips([12]). In this paper, we give
$p$-adic $q$-integral representation for Kim's $q$-Bernstein
polynomials and investigate some interesting identities of
$q$-Bernstein polynomials associated with $q$-extension of binomial
distribution, $q$-Stirling numbers and Carlitz's $q$-Bernoulli
numbers.}

\section{Introduction}

Let $p$ be a fixed prime number. Throughout this paper, $\Bbb Z_p$,
$\Bbb Q_p$, $\Bbb C$ and $\Bbb C_p$ denote the ring of $p$-adic
integers, the field of $p$-adic rational numbers, the complex number
field and the completion of algebraic closure of $\Bbb Q_p$,
respectively. Let $\Bbb N$ be the set of natural numbers and $\Bbb
Z_+ =\Bbb N \cup \{0\}$. Let $\nu_p$ be the normalized exponential
valuation of $\Bbb C_p$ with $|p|_p=p^{-\nu_p (p)} =\frac{1}{p}$.

When one talks of $q$-extension, $q$ is variously considered as an
indeterminate, a complex number $q \in \Bbb C$ or $p$-adic number $q
\in \Bbb C_p$. If $q\in \Bbb C$, one normally assumes $|q|<1$, and
if $q\in \Bbb C_p$, one normally assumes $|1-q|_p <1$. 

The $q$-bosonic natural numbers are defined by
$[n]_q=\frac{1-q^n}{1-q} =1+q+q^2+\cdots+q^{n-1}$ for $n \in \Bbb
N$, and the $q$-factorial is defined by $[n]_q !=[n]_q [n-1]_q
\cdots [2]_q [1]_q$. For the $q$-extension of binomial coefficient,
we use the following notation in the form of
$$\binom{n}{k}_q=\frac{[n]_q!}{[n-k]_q![k]_q!}=\frac{[n]_q [n-1]_q
\cdots [n-k+1]_q}{[k]_q!}.
$$

Let $C[0, 1]$ denote the set of continuous functions on $[0,
1](\subset \mathbb{R})$. Then Bernstein operator for $f \in C[0,1]$
is defined by
\begin{eqnarray*}
\mathbb{B}_n(f | x)= \sum_{k=0}^{n} f(\frac{k}{n})\binom{n}{k} x^k
(1-x)^{n-k}=\sum_{k=0}^{n} f(\frac{k}{n})B_{k,n}(x),
\end{eqnarray*}
where $n, k \in \Bbb Z_+$. The polynomials $B_{k,n}(x)=
\binom{n}{k}x^k (1-x)^{n-k}$ are called Bernstein polynomials of
degree $n$ (see [1]). For $f \in C[0, 1]$, Kim's $q$-Bernstein
operator of order $n$ for $f$ is defined by
\begin{eqnarray*}
\mathbb{B}_{n,q} (f | x)= \sum_{k=0}^{n} f(\frac{k}{n})\binom{n}{k}
[x]_q^k [1-x]_{\frac{1}{q}}^{n-k}=\sum_{k=0}^{n}
f(\frac{k}{n})B_{k,n}(x,q),
\end{eqnarray*}
where $n, k \in \Bbb Z_+$. Here $B_{k,n}(x,q)= \binom{n}{k}[x]_q^k
[1-x]_{\frac{1}{q}}^{n-k}$ are called the Kim's $q$-Bernstein polynomials
of degree $n$ (see [4]).

We say that $f$ is uniformly differentiable function at a point $a
\in \Bbb Z_p$, and write $f\in UD(\Bbb Z_p)$, if the difference
quotient $F_f (x,y) =\frac{f(x)-f(y)}{x-y}$ has a limit $f'(a)$ as
$(x,y)\rightarrow (a,a)$. For $f\in UD(\Bbb Z_p)$, the $p$-adic
$q$-integral on $\Bbb Z_p$ is defined by
\begin{eqnarray*}
I_{q} (f)=\int_{\Bbb Z_p } f(x) d\mu_{q} (x) = \lim_{N \rightarrow
\infty} \frac{1}{[p^N]_q} \sum_{x=0}^{p^N-1} f(x)q^x, \quad
(\text{see [6]}).
\end{eqnarray*}

Carlitz's $q$-Bernoulli numbers can be represented by $p$-adic
$q$-integral on $\Bbb Z_p$ as follows:
\begin{eqnarray}
\int_{\Bbb Z_p } [x]_q^n d\mu_{q} (x) = \lim_{N \rightarrow \infty}
\frac{1}{[p^N]_q} \sum_{x=0}^{p^N-1} [x]_q^n q^x = \beta_{n, q},
\quad (\text{see [6, 7]}).
\end{eqnarray}
The $k$-th order factorial of the $q$-number $[x]_q$, which is
defined by $$[x]_{k,q}=[x]_q [x-1]_q \cdots [x-k+1]_q =
\frac{(1-q^x)(1-q^{x-1})\cdots(1-q^{x-k+1})}{(1-q)^k},$$ is called
the $q$-factorial of $x$ of order $k$ (see [6]).

In this paper, we give $p$-adic $q$-integral representation for
Kim's $q$-Bernstein polynomials and derive some interesting
identities for the Kim's $q$-Bernstein polynomials associated with
$q$-extension of binomial distribution, $q$-Stirling numbers and
Carlitz's $q$-Bernoulli numbers.

\medskip

\section{$q$-Bernstein polynomials}

In this section, we assume that $0<q<1$. Let $\mathbb{P}_q = \{
{\underset{i}{\sum}}a_i [x]_q^i | \,\, a_i \in \mathbb{R} \}$ be the
space of $q$-polynomials of degree less than or equal to $n$.

For $f \in C[0, 1]$ and $n, k \in \Bbb Z_+$, Kim's $q$-Bernstein
operator of order $n$ for $f$ is defined by
\begin{eqnarray}
\mathbb{B}_{n,q} (f | x)=\sum_{k=0}^{n} f(\frac{k}{n})B_{k,n}(x,q).
\end{eqnarray}
Here $B_{k,n}(x,q)= \binom{n}{k}[x]_q^k [1-x]_{\frac{1}{q}}^{n-k}$
are the Kim's $q$-Bernstein polynomials of degree $n$ (see [4]).

Kim's $q$-Bernstein polynomials of degree $n$ is a basis for the
space of $q$-polynomials of degree less than or equal to $n$. That
is, Kim's $q$-Bernstein polynomials of degree $n$ is a basis for
$\mathbb{P}_q$.

We see that Kim's $q$-Bernstein polynomials of degree $n$
span the space of $q$-polynomials. That is, any $q$-polynomials of
degree less than or equal to $n$ can be written as a linear
combination of the Kim's $q$-Bernstein polynomials of degree $n$.
For $n, k \in \Bbb Z_+$ and $x \in [0, 1]$, we have
\begin{eqnarray}B_{k,n}(x,q)=\sum_{l=k}^{n} \binom{n}{l}\binom{l}{k}(-1)^{l-k}[x]_q^l, \quad (\text{see [4]}).\end{eqnarray}

If there exist constants $C_0, C_1, \ldots, C_n$ such that $C_0
B_{0,n}(x,q)+C_1 B_{1,n}(x,q)+\cdots+C_n B_{n,n}(x,q)=0$ holds for
all $x$, then we can derive the following equation from (3):
\begin{eqnarray*}
0&=&C_0 B_{0,n}(x,q)+C_1 B_{1,n}(x,q)+\cdots+C_n B_{n,n}(x,q)\\
&=&C_0 \sum_{i=0}^{n}(-1)^i \binom{n}{i}\binom{i}{0}[x]_q^i+C_1
\sum_{i=1}^{n}(-1)^{i-1} \binom{n}{i}\binom{i}{1}[x]_q^i \\ & &
\qquad \qquad \qquad \qquad \qquad \qquad \qquad +\cdots+
C_n \sum_{i=n}^{n}(-1)^{i-n} \binom{n}{i}\binom{i}{n}[x]_q^i\\
&=&C_0 + \{ \sum_{i=0}^{1}C_i (-1)^{i-1} \binom{n}{1}\binom{1}{i}
\}[x]_q +\cdots+ \{ \sum_{i=0}^{n}C_i (-1)^{i-n}
\binom{n}{n}\binom{n}{i} \}[x]_q^n.
\end{eqnarray*}
Since the power basis is a linearly independent set, it follows that
\begin{eqnarray*}
C_0  \qquad \qquad=0,\\
\sum_{i=0}^{1}C_i (-1)^{i-1} \binom{n}{1}\binom{1}{i}=0,\\
\vdots \qquad \qquad \qquad \vdots \,\, \\ \sum_{i=0}^{n}C_i
(-1)^{i-n} \binom{n}{n}\binom{n}{i}=0,
\end{eqnarray*}
which implies that $C_0=C_1=\cdots=C_n =0$ ($C_0$ is clearly zero,
substituting this in the second equation gives $C_1 =0$,
substituting these two into the third equation  gives $C_2=0$, and
so on).

Let us consider a $q$-polynomial $P_q (x) \in \mathbb{P}_q$ which is
written by a linear combination of Kim's $q$-Bernstein basis
functions as follows:
\begin{eqnarray}
P_q (x)=C_0 B_{0,n}(x,q)+C_1 B_{1,n}(x,q)+\cdots+C_n B_{n,n}(x,q).
\end{eqnarray}
It is easy to write (4) as a dot product of two values.
\begin{eqnarray}
P_q (x)=(B_{0,n}(x,q), B_{1,n}(x,q), \ldots, B_{n,n}(x,q)) \left(
\begin{array}{llll}C_0 \\C_1 \\ \, \vdots \\ C_n \end{array} \right).
\end{eqnarray}
From (5), we can derive the following equation:
\begin{eqnarray*}
P_q (x)=(1, [x]_q, \ldots, [x]_q^n) \left(
\begin{array}{lllll}b_{00} \ \  0 \ \ \ 0 \ \ \cdots \ \ 0\\b_{10} \ \ b_{11} \ 0  \ \  \cdots \ \ 0
\\ \ \ \vdots \ \ \ \ \vdots \ \ \ \ \vdots \ \  \ddots \ \ \vdots \\b_{n0} \  b_{n1} \ b_{n2} \
\ldots  b_{nn} \end{array} \right) \left(
\begin{array}{llll}C_0 \\C_1 \\ \, \vdots \\ C_n \end{array}
 \right),
\end{eqnarray*}
where the $b_{ij}$ are the coefficients of the power basis that are
used to determine the respective Kim's $q$-Bernstein polynomials. We
note that the matrix in this case is lower triangular.

From (2) and (3),  we note that
\begin{eqnarray*}
B_{0, 2}(x,q)&=&[1-x]_{\frac{1}{q}}^2=\sum_{l=0}^{2}\binom{2}{l}(-1)^l [x]_q^l=1-2[x]_q+[x]_q^2, \\
B_{1, 2}(x,q)&=&\binom{2}{1}[x]_q [1-x]_{\frac{1}{q}}=2[x]_q-2[x]_q^2, \\
B_{2, 2}(x,q)&=&\binom{2}{2}[x]_q^2=[x]_q^2. \\
\end{eqnarray*}

In the quadratic case ($n=2$), the matrix representation is
\begin{eqnarray*}
P_q (x)=(1, \, [x]_q, \, [x]_q^2) \left(
\begin{array}{lll} \ \ 1 \ \ \ \ \, 0 \ \ \ \, 0 \\-2 \ \ \ \  2 \ \ \ \, 0
\\\ \ 1 \   -2 \ \ \ 1 \end{array} \right) \left(
\begin{array}{llll}C_0 \\C_1 \\ C_2 \end{array}
 \right).
\end{eqnarray*}
In the cubic case ($n=3$), the matrix representation is
\begin{eqnarray*}
P_q (x)=(1, \, [x]_q, \, [x]_q^2, \, [x]_q^3) \left(
\begin{array}{llll} \ \ 1 \ \ \ \ \, 0 \ \ \ \ \, 0 \ \ \  0 \\-3 \ \ \ \ 3 \ \ \
\ \,0  \ \ \   0 \\\ \ 3 \   -6 \ \ \ \ 3 \  \ \  0 \\-1 \ \ \ \ 3 \
-3 \ \ \, 1 \end{array} \right) \left(
\begin{array}{llll}C_0 \\C_1 \\ C_2 \\ C_3\end{array}
 \right).
\end{eqnarray*}
In many applications of $q$-Bernstein polynomials, a matrix
formulation for the Kim's $q$-Bernstein polynomials seems to be
useful.

\medskip

\section{$q$-Bernstein polynomials, $q$-Stirling numbers and $q$-Bernoulli numbers}

In this section, we assume that $q \in \Bbb C_p$ with $|1-q|_p<1$.

For $f \in UD(\Bbb Z_p)$, let us consider the $p$-adic analogue of
Kim's $q$-Bernstein type operator of order $n$ on $\Bbb Z_p$ as
follows:
\begin{eqnarray*}
\mathbb{B}_{n,q} (f | x)= \sum_{k=0}^{n} f(\frac{k}{n})\binom{n}{k}
[x]_q^k [1-x]_{\frac{1}{q}}^{n-k}=\sum_{k=0}^{n}
f(\frac{k}{n})B_{k,n}(x,q).
\end{eqnarray*}

Let $(Eh)(x)=h(x+1)$ be the shift operator. Then the $q$-difference
operator is defined by
\begin{eqnarray}
\Delta_q^n:=(E-I)_q^n = \prod_{i=1}^n (E-q^{i-1}I),
\end{eqnarray}
where $(Ih)(x)=h(x)$. From (6), we derive the following equation:
\begin{eqnarray}
\Delta_q^n f(0)=\sum_{k=0}^{n} {\binom{n}{k}}_q (-1)^k
q^{\binom{k}{2}}f(n-k), \quad (\text{see [7]}).
\end{eqnarray}
By (7), we easily see that
$$f(x)=\sum_{n \ge 0} \binom{x}{n}_q \Delta_q^n f(0), \quad (\text{see [6, 7]}).$$

The $q$-Stirling number of the first kind is defined by
\begin{eqnarray}
\prod_{k=1}^n (1+[k]_q z)= \sum_{k=0}^{n} S_{1,q}(n,k)z^k, \quad
(\text{see [5, 6]}),
\end{eqnarray}
and the $q$-Stirling number of the second kind is also defined by
\begin{eqnarray}
\prod_{k=1}^n (\frac{1}{1+[k]_q z})= \sum_{k=0}^{n} S_{2,q}(n,k)z^k,
\quad (\text{see [5]}).
\end{eqnarray}
By (6), (7), (8) and (9), we get
\begin{eqnarray*}
S_{2,q}(n,k)&=& \frac{q^{- \binom{k}{2}}}{[k]_q !} \sum_{j=0}^{k}
(-1)^j
q^{\binom{j}{2}}\binom{k}{j}_q [k-j]_q^n \\
&=& \frac{q^{- \binom{k}{2}}}{[k]_q !} \Delta_q^k 0^n,
\end{eqnarray*}
for $n, k \in \Bbb Z_+$ (see [6]).

Let us consider Kim's $q$-Bernstein polynomials of degree $n$ on
$\Bbb Z_p$ as follows: $$B_{k,n}(x,q)= \binom{n}{k}[x]_q^k
[1-x]_{\frac{1}{q}}^{n-k},$$ for $n, k \in \Bbb Z_+$ and $x \in \Bbb
Z_p$. Thus, we easily see that
\begin{eqnarray}
\int_{\Bbb Z_p} B_{k,n}(x, q) d\mu_q(x) = \sum_{l=0}^{n-k}
\binom{n-k}{l} \binom{n}{k} (-1)^{l} \int_{\Bbb Z_p}[x]_q^{l+k} d
\mu_q(x).
\end{eqnarray}
By (1) and (10), we obtain the following proposition.

\begin{proposition}
For $n, k \in \Bbb Z_+$, we have
\begin{eqnarray*}
\int_{\Bbb Z_p} B_{k,n}(x, q) d\mu_q(x) =  \sum_{l=0}^{n-k}
\binom{n-k}{l} \binom{n}{k} (-1)^{l} \beta_{l+k, q},
\end{eqnarray*}
where $\beta_{l+k, q}$ are the $(l+k)$-th Carlitz's $q$-Bernoulli
numbers.
\end{proposition}

\medskip

From the definition of Kim's $q$-Bernstein polynomial, we note that
\begin{eqnarray}
\sum_{k=i}^{n} \frac{\binom{k}{i}}{\binom{n}{i}} B_{k,n}(x, q) =
\sum_{k=0}^{i} q^{\binom{k}{2}}\binom{x}{k}_q[k]_q! S_{2,q}(k, i-k),
\end{eqnarray}
where $i \in \Bbb N$. From the definition of $q$-binomial
coefficient, we have
\begin{eqnarray}
\binom{n+1}{k}_q = \binom{n}{k-1}_q + q^k \binom{n}{k}_q =q^{n-k}
\binom{n}{k-1}_q + \binom{n}{k}_q.
\end{eqnarray}
By (12), we see that
\begin{eqnarray}
\int_{\Bbb Z_p} \binom{x}{n}_q d \mu_q (x) =
\frac{(-1)^n}{[n+1]_q}q^{(n+1)-\binom{n+1}{2}}, \quad \text{(see [6,
7])}.
\end{eqnarray}
From (1), (11) and (13), we obtain the following theorem.

\begin{theorem}
For $n, k \in \Bbb Z_+$ and $i \in \Bbb N$, we have
\begin{eqnarray*}
& & \sum_{k=i}^{n}\sum_{l=0}^{n-k}
\frac{\binom{k}{i}}{\binom{n}{i}}\binom{n-k}{l} \binom{n}{k}(-1)^l
\beta_{l+k, q} \\& & \qquad = \sum_{k=0}^{i} q^{\binom{k}{2}}[k]_q!
S_{2,q}(k, i-k)\frac{(-1)^k}{[k+1]_q}q^{(k+1)-\binom{k+1}{2}}.
\end{eqnarray*}
\end{theorem}

\medskip

It is easy to see that for $i \in \Bbb N$,
\begin{eqnarray}
& & \sum_{k=i}^{n} \frac{\binom{k}{i}}{\binom{n}{i}}B_{k,n}(x,q)
=[x]_q^i.
\end{eqnarray}
By (11) and (14), we easily get
\begin{eqnarray*}
[x]_q^i = \sum_{k=0}^{i} q^{\binom{k}{2}} \binom{x}{k}_q [k]_q!
S_{2,q}(k, i-k), \quad \text{(see [6])}.
\end{eqnarray*}
Thus, we have
\begin{eqnarray}
& & \int_{\Bbb Z_p}[x]_q^i d \mu_q (x)= \sum_{k=0}^{i}
q^{\binom{k}{2}}
[k]_q! S_{2,q}(k, i-k) \int_{\Bbb Z_p} \binom{x}{k}_q d \mu_q (x) \\
& & \qquad =q \sum_{k=0}^i [k]_q!S_{2,q}(k,
i-k)\frac{(-1)^k}{[k+1]_q}. \notag
\end{eqnarray}
By (1) and (15), we obtain the following corollary.

\begin{corollary}
For $n, k \in \Bbb Z_+$ and $i \in \Bbb N$, we have
\begin{eqnarray*}
\beta_{i, q}=q \sum_{k=0}^i [k]_q!S_{2,q}(k,
i-k)\frac{(-1)^k}{[k+1]_q}.
\end{eqnarray*}
\end{corollary}

\medskip

It is known that
\begin{eqnarray}
S_{2,q}(n,k)=\frac{1}{(1-q)^k}\sum_{j=0}^k (-1)^{k-j}
\binom{k+n}{k-j}\binom{j+n}{j}_q, \quad \text{(see [6])},
\end{eqnarray}
and $$\binom{n}{k}_q = \sum_{j=0}^{n}
\binom{n}{j}(q-1)^{j-k}S_{2,q}(k, j-k).
$$
By simple calculation, we have that
\begin{eqnarray}
q^{nx}&=& \sum_{k=0}^{n} (q-1)^k q^{\binom{k}{2}} \binom{n}{k}_q
[x]_{k,q}\\
&=& \sum_{m=0}^{n} \{ \sum_{k=m}^{n}(q-1)^k \binom{n}{k}_q
S_{1,q}(k,m) \} [x]_q^m \notag
\end{eqnarray}
and
\begin{eqnarray}q^{nx} =\sum_{m=0}^{n} \binom{n}{m} (q-1)^m
[x]_q^m.\end{eqnarray}
From (17) and (18), we note that
\begin{eqnarray*}
\binom{n}{m} =\sum_{k=m}^{n}(q-1)^{-m+k} \binom{n}{k}_q
S_{1,q}(k,m), \quad (\text{see [6]}).
\end{eqnarray*}
Thus, we obtain the following proposition.

\begin{proposition}
For $n, k \in \Bbb Z_+$, we have $$B_{k,n}(x,q)= \binom{n}{k}[x]_q^k
[1-x]_{\frac{1}{q}}^{n-k}=\sum_{m=k}^{n}(q-1)^{-k+m} \binom{n}{m}_q
S_{1,q}(m, k)[x]_q^k[1-x]_{\frac{1}{q}}^{n-k}.$$
\end{proposition}

\medskip

From the definition of the $q$-Stirling numbers of the first kind,
we get
\begin{eqnarray}
q^{\binom{n}{2}} \binom{x}{n}_q [n]_q!=[x]_{n,q} \,
q^{\binom{n}{2}}=\sum_{k=0}^n S_{1,q}(n, k)[x]_q^k.
\end{eqnarray}
By (11) and (19), we obtain the following theorem.

\begin{theorem}
For $n, k \in \Bbb Z_+$ and $i \in \Bbb N$, we have
\begin{eqnarray*}
& & \sum_{k=i}^{n} \frac{\binom{k}{i}}{\binom{n}{i}}B_{k,n}(x,q)=
 \sum_{k=0}^{i} \sum_{l=0}^{k} S_{1,q}(k, l)S_{2,q}(k, i-k)[x]_q^l.
\end{eqnarray*}
\end{theorem}

\medskip

By (14) and Theorem 5, we obtain the following corollary.

\begin{corollary}
For $i \in \Bbb Z_+$, we have
\begin{eqnarray*}
\beta_{i, q}=\sum_{k=0}^{i} \sum_{l=0}^{k} S_{1,q}(k, l)S_{2,q}(k,
i-k)\beta_{l.q}.
\end{eqnarray*}
\end{corollary}

The $q$-Bernoulli polynomials of order $k \in \Bbb Z_+$ are defined
by
\begin{eqnarray}
\beta_{n,q}^{(k)}(x)=\frac{1}{(1-q)^n}\sum_{i=0}^n
\binom{n}{i}(-1)^i q^{ix} \int_{\Bbb Z_p}\cdots  \int_{\Bbb Z_p}
q^{\sum_{l=1}^k (k-l+i)x_l} d\mu_q(x_1) \cdots d\mu_q(x_k).
\end{eqnarray}
Thus, we have
\begin{eqnarray*}
\beta_{n,q}^{(k)}(x)=\frac{1}{(1-q)^n}\sum_{i=0}^n \frac{
(-1)^i\binom{n}{i}(i+k)\cdots (i+1)}{[i+k]_q \cdots [i+1]_q} q^{ix},
\quad (\text{see [6]}).
\end{eqnarray*}

The inverse $q$-Bernoulli polynomials of order $k$ are defined by
\begin{eqnarray}
\beta_{n,q}^{(-k)}(x)=\frac{1}{(1-q)^n}\sum_{i=0}^n \frac{(-1)^i
\binom{n}{i}q^{ix}}{\int_{\Bbb Z_p}\cdots  \int_{\Bbb Z_p}
q^{\sum_{l=1}^k (k-l+i)x_l} d\mu_q(x_1) \cdots d\mu_q(x_k)}, \quad
(\text{see [6]}).
\end{eqnarray}
In the special case $x=0$, $\beta_{n,q}^{(k)}(0)= \beta_{n,q}^{(k)}$
are called the $n$-th $q$-Bernoulli numbers of order $k$ and
$\beta_{n,q}^{(-k)}(0)= \beta_{n,q}^{(-k)}$ are also called the
inverse $q$-Bernoulli numbers of order $k$.

From (21), we have
\begin{eqnarray}
\beta_{k, q}^{(-n)} &=& \frac{1}{(1-q)^k}\sum_{j=0}^k (-1)^j
\binom{k}{j} \frac{[j+n]_q \cdots [j+1]_q}{(j+n)\cdots(j+1)} \notag \\
&=& \frac{1}{(1-q)^k}\sum_{j=0}^k (-1)^j
\frac{\binom{k+n}{n-j}}{\binom{k+n}{n}}\binom{j+n}{n}_q \frac{[n]_q
!}{n!}\\
&=& \frac{[n]_q!}{\binom{k+n}{n}n!} \{ \frac{1}{(1-q)^k}\sum_{j=0}^k
(-1)^j \binom{k+n}{n-j}\binom{j+n}{n}_q \}. \notag
\end{eqnarray}
By (16) and (22), we get
\begin{eqnarray}
\frac{n!}{[n]_q!}\binom{k+n}{n} \beta_{k,q}^{(-n)} = S_{2, q}(n, k).
\end{eqnarray}
Therefore, by (11) and (23), we obtain the following theorem.

\medskip

\begin{theorem}
For $i, n, k \in \Bbb Z_+$, we have
\begin{eqnarray*}
& & \sum_{k=i}^{n} \frac{\binom{k}{i}}{\binom{n}{i}}B_{k,n}(x,q)=
 \sum_{k=0}^{i} q^{\binom{k}{2}}k! \binom{i}{k}\binom{x}{k}_q
 \beta_{i-k, q}^{(-k)}.
\end{eqnarray*}
\end{theorem}

\medskip

It is easy to show that
\begin{eqnarray*}
q^{\binom{n}{2}}\binom{x}{n}_q &=& \frac{1}{[n]_q!}
\prod_{k=0}^{n-1}([x]_q-[k]_q) \\&=&
\frac{1}{[n]_q!}\sum_{k=0}^{n}(-1)^k [x]_q^{n-k} S_{1,q}(n-1, k).
\end{eqnarray*}
Thus, we have that
\begin{eqnarray*}
& & \sum_{k=i}^{n} \frac{\binom{k}{i}}{\binom{n}{i}}B_{k,n}(x,q)=
 \sum_{k=0}^{i} \sum_{j=0}^{k} (-1)^j [x]_q^{k-j} S_{1,q}(k-1, j)
 \frac{k!}{[k]_q!}  \binom{i}{k} \beta_{i-k, q}^{(-k)},
\end{eqnarray*}
where $n, k, i \in \Bbb Z_+$.

\medskip

\noindent \textbf{Acknowledgement.}
This paper was supported by the research grant of Kwangwoon University in 2010.

\end{document}